\documentclass[11pt]{amsart}

\usepackage{amsmath}
\usepackage{amsfonts}
\usepackage{amssymb}
\usepackage{latexsym}
\usepackage[dvips]{graphicx}

\newtheorem{Thm}{Theorem}[section]

\newtheorem{Prop}[Thm]{Proposition}

\newtheorem{Lem}[Thm]{Lemma}
\newtheorem{Def}[Thm]{Definition}

\def\bR {\mathbf{R}}

\def\cD {\mathcal{D}}

\def\a {{\alpha}}
\def\b {{\beta}}

\def\g {{\gamma}}

\def\eps {{\epsilon}}

\def\th {{\theta}}

\def\vphi{{\varphi}}

\def\Om {{\Omega}}

\def\rstr {{\big |}}
\def\indc {{\bf 1}}

\def\d {{\partial}}

\def\wto {\rightharpoonup}

\newcommand{\Div}{\operatorname{div}}

\newcommand{\Curl}{\operatorname{curl}}

\newcommand{\be}{\begin{equation}}
\newcommand{\ee}{\end{equation}}
\newcommand{\ba}{\begin{aligned}}
\newcommand{\ea}{\end{aligned}}
\newcommand{\lb}{\label}

\begin{document}

\title[Diffusion limit for generalized Carleman]
{The nonlinear diffusion limit\\for generalized Carleman models: \\
the initial-boundary value problem}

\bibliographystyle{plain}
\author[F. Golse]{Fran\c cois Golse}
\address{F.G.: Universit\'e Paris 7 - Denis Diderot, Laboratoire J.-L. Lions, 
Bo\^\i te courrier 187, 4 place Jussieu, 
75252 Paris Cedex 05, France} 

\author[F. Salvarani]{Francesco Salvarani}
\address{F.S.: Dipartimento di Matematica,
Universit\`a degli Studi di Pavia.
Via Ferrata 1, 27100 Pavia, Italy}

\begin{abstract}
Consider the initial-boundary value problem for the 2-speed Carleman model of the
Boltzmann equation of the kinetic theory of gases (see [T. Carleman, ``Probl\` emes  
math\' ematiques dans la th\' eorie cin\' etique des gaz", Almqvist-Wiksells, Uppsala, 
1957]), set in some bounded interval with
boundary conditions prescribing the density of particles entering the interval. Under
the usual parabolic scaling, a nonlinear diffusion limit is established for this problem.
In fact, the techniques presented here allow treating generalizations of the Carleman 
system where the collision frequency is proportional to the $\a$-th power of the 
macroscopic density, with $\a\in[-1,1]$.
\end{abstract}

\maketitle

\section{Carleman models and their diffusion limits}

In the 1930's, Carleman proposed a model \cite{Carl} describing the time evolution of 
a monodimensional gas composed of two species of particles that move at a constant
speed $c>0$ in the $x$-direction.  The number density at time $t$ and position $x$ of 
particles moving at speed $+c$ is denoted by $u=u(x,t)$ while that of particles moving 
at speed $-c$ is denoted $v=v(x,t)$. Carleman's system is
\be\lb{OrigCarl}
\ba
\d_tu + c\d_xu=(u+v)(v-u)\,,
\\
\d_tv - c\d_xv=(u+v)(u-v)\,.
\ea
\ee
A very natural extension of the model, called the \textit{generalized Carleman model},
involves a collision frequency that is proportional to some power of the macroscopic 
density $\rho=u+v$, as follows:
\be\lb{GenCarl}
\ba
\d_tu + c\d_xu=(u+v)^\a(v-u)\,,
\\
\d_tv - c\d_xv=(u+v)^\a(u-v)\,.
\ea
\ee
For $\a=1$ the original Carleman system is recovered, whereas $\a=0$ gives another 
remarkable system, known as the Goldstein-Taylor model \cite{Gol,Tay}, that can be
reduced to a damped wave equation (the telegrapher's equation). In the latter case,
there is an explicit representation of the solution in terms of the standard Poisson: see
\cite{McKean}, p. 56 for this formula, originally found by M. Kac. 

Other instances of kinetic models with a collision frequency whose dependence on 
the macroscopic density is other than linear can be found in Radiative Transfer: see
for instance \cite{Mihalas}.

The only nontrivial hydrodynamic limits of the Carleman model, generalized or not, 
are diffusion limits --- linear or nonlinear. Indeed, local equilibria for those models
are those densities for which $(u-v)(u+v)^\alpha=0$, implying $u=v$. Hence all local 
equilibria for those models have mean velocity $0$: in the language of the kinetic
theory of gases, analogues for the system (\ref{OrigCarl}) of local Maxwellians have 
mean velocity $0$. Hence the limiting dynamics of the system can be observed only 
on a longer time scale that corresponds with a diffusion equation.

After setting $c=1$ without loss of generality, we consider the limit as $\eps\to 0^+$ of 
the following scaled version of (\ref{GenCarl}):
\be\lb{ScldCrlm}
\ba
\eps^2\d_tu_\eps + \eps\d_xu_\eps=(u_\eps+v_\eps)^\a(v_\eps-u_\eps)\,,
\\
\eps^2\d_tv_\eps - \eps\d_xv_\eps=(u_\eps+v_\eps)^\a(u_\eps-v_\eps)\,.
\ea
\ee
Defining the macroscopic mass density $\rho_\eps$ and the current $j_\eps$ by
\be
\rho_\eps = u_\eps+ v_\eps\,, \qquad  j_\eps= \frac{u_\eps -v_\eps}{\eps}\,,
\ee
we put the system (\ref{ScldCrlm}) in the form
\be\lb{MacroScldCarl}
\ba
\d_t\rho_\eps + \d_xj_\eps&=0\,,
\\
\eps^2\d_tj_\eps + \d_x\rho_\eps&=-2\rho_\eps^\a j_\eps\,.
\ea
\ee
Then, we prove that the term $\eps^2\d_tj_\eps$ is negligeable in the vanishing $\eps$
limit, and show that the limiting density is governed by the following nonlinear diffusion
equation
\be
\label{NLDiff}
\d_t\rho=\tfrac12\d_{xx}\left(\frac{\rho^{1-\a}}{1-\a}\right)\,,
\ee
for $\a\in[-1,1)$, while the case $\a=1$ leads to
\be
\label{NLDiff-1}
\d_t\rho=\tfrac12\d_{xx}\ln\rho\,.
\ee

Several result on this problem have been obtained over the last thirty years.

In particular, existence and uniqueness for the solution of the initial-boundary value problem
of (\ref{OrigCarl}) have been proven by Fitzgibbon \cite{Fitzg} --- see also \cite{Kolod} for the 
same problem on the infinite line, and \cite{ReedIllner1}, \cite{ReedIllner2} and \cite{Tartar} for 
more information concerning the large time behavior of the solutions\footnote{We are grateful
to L. Tartar for indicating the references \cite{Kolod} and \cite{ReedIllner1}, \cite{ReedIllner2}.}.

Several authors (among others Kurtz \cite{Kurtz}, McKean \cite{McKean}, Fitzgibbon \cite{Fitzg},
Pulvirenti and Toscani \cite{PuTo}, P.-L. Lions and Toscani \cite{LionsToscani}, Marcati and 
Rubino \cite{MaRu}, Donatelli and Marcati \cite{DoMa}, Salvarani and V\'azquez \cite{SaVa}) 
have considered the relaxation problem for the Carleman system or generalizations thereof.

All these papers on the diffusion limit for kinetic models deal with the initial value problem on
the infinite line, or with the initial-boundary value problem with specular or periodic conditions
at the boundary.

Establishing the diffusion limit for Carleman type models (generalized or not) in a bounded 
domain with non homogeneous boundary conditions is more difficult. Indeed, mimicking the 
proofs of the diffusion limit in the cases mentioned above fails to produce a uniform bound 
on the current $j_\eps$ in the case of a nonhomogeneous boundary value problem.

To the best of our knowledge, the only available convergence proof for an initial-boundary
value problem with very general boundary conditions at the time of this writing is in \cite{Salv} 
for the Goldstein-Taylor model. Unfortunately, the proof in \cite{Salv} uses extensively the 
linear nature of the problem coming from the assumption $\a=0$, and does not seem to be
extendible to nonlinear models.

The present work establishes the diffusion limit for all Carleman systems (\ref{GenCarl})
with $\a\in[-1,1]$, in a bounded domain $\Om=(0,1)$, with boundary conditions imposing
the density of particles entering the domain $\Om$
\be\lb{BoundCond}
\ba
u_\eps(t,0)=\vphi^-(t)\,,\qquad t>0\,, 
\\
v_\eps(t,1)=\vphi^+(t)\,,\qquad t>0\,.
\ea
\ee
We have restricted our analysis to the case $\vert\alpha\vert\leq 1$, as it leads to a unified 
treatment based on the dissipative nature of the problem.

The new ingredient in the present paper is a uniform estimate bearing on some notion of
relative entropy of the solution $(u_\eps,v_\eps)$ with respect to a well-chosen profile that
satisfies the boundary conditions expected to hold in the vanishing $\eps$ limit. By failing
to introducing this profile, previous attempts to establishing the nonlinear diffusion limit in
the case of nonhomogeneous data fell short of obtaining uniform bounds on the current 
$j_\eps$, except in the particular case $\varphi^+(t)=\varphi^-(t)= {\it Const}$.

This idea of using the relative entropy with respect to some adequate profile in order to
match nontrivial boundary conditions can be used on more complicated models. In the
introduction above, we already mentioned that the generalized Carleman equations are
somewhat analogous to more complicated models appearing in Radiative Transfer, for
which the nonlinear diffusion approximation is known under the name of ``Rosseland
approximation". We shall apply the method presented here to these more complicated
models in a forthcoming paper \cite{GolSal}.

\section{Main results}

Before stating the convergence theorem that is the main result in this paper, we recall
some background on the Carleman systems.

\subsection{Existence and uniqueness theory for Carleman systems}

The class of initial and boundary data for the Carleman systems considered in this paper,
which we henceforth call ``admissible data", is defined below.

\begin{Def}\lb{DEF-ADMIS}
Consider the system $(\ref{ScldCrlm})$ with $|\a|\le 1$ posed for $(t,x)\in(0,T)\times(0,1)$,
with boundary conditions $(\ref{BoundCond})$ and initial condition
\be\lb{InCond}
\ba
u_\eps(0,x)=u^{in}(x)\,,\qquad 0<x<1\,, 
\\
v_\eps(0,x)=v^{in}(x)\,,\qquad 0<x<1\,.
\ea
\ee
The initial data $(u^{in},v^{in})$ and the boundary data $(\vphi^-,\vphi^+)$ are said to be
admissible on the time interval $(0,T)$ if and only if
\begin{enumerate}
\item $\vphi^\pm\in W^{1,\infty}([0,T])$ and $\vphi^\pm>0$  on $(0,T)$, while
\item $u^{in}$ and $v^{in}\in L^3(0,1)$ and $u^{in},v^{in}\ge 0$ a.e. in $(0,1)$.
\end{enumerate}
\end{Def}

The existence and uniqueness of a nonnegative solution for the original Carleman
model is well-known. The following result is a standard generalization of Fitzgibbon's
in \cite{Fitzg2}:

\begin{Thm}
Consider the generalized Carleman model $(\ref{GenCarl})$ for $|\a|\le 1$ in the domain
$\Om= (0,1)$ with boundary condition $(\ref{BoundCond})$ and initial condition $(\ref{InCond})$.
Assume that the initial and boundary data $(u^{in},v^{in})$ and $\vphi^\pm$ are admissible.
Then there exists a unique nonnegative generalized solution $(u,v)$ of $(\ref{GenCarl})$
in $C([0,T];L^1(0,1)\times L^1(0,1))$.
\end{Thm}

\begin{proof}
It is well known that the unbounded operator
$$
B_\a(u,v)=\left(-\d_xu+(u+v)^\a(v-u),\d_xv+(u+v)^\a(u-v)\right)
$$
on $L^{1}(0,1)\times L^{1}(0,1)$ with domain
$$
\cD(B_\a)=\{(u,v)\in W^{1,1}(0,1)\times W^{1,1}(0,1)\,|\,u(0)=0\hbox{ and }v(1)=0\}
$$
is dissipative (see, for example, \cite{LionsToscani, SaVa}). The existence and uniqueness of the
solution $(u,v)$ easily follows from the same method as in \cite{Fitzg2}.
\end{proof}

\subsection{Uniqueness theory for the nonlinear diffusion equation}

Next we consider the target, nonlinear diffusion equation (\ref{NLDiff}) (or (\ref{NLDiff-1})
when $\a=1$).

\begin{Def}\lb{DEF-ADMIS-DIFF-CL}
Consider the nonlinear diffusion equation $(\ref{NLDiff})$ or $(\ref{NLDiff-1})$ when $\a=1$
with Dirichlet boundary conditions
\be\lb{BCDiff}
\rho(t,0)=2\varphi^-\,,\quad\rho(t,1)=2\varphi^+\,,\qquad t\in(0,T)\,,
\ee
and initial condition
\be\lb{ICDiff}
\rho(0,x)=\rho_0(x)\,,\qquad x\in(0,1)\,.
\ee
We call these data admissible if and only if 
\begin{enumerate}
\item for $\a\in[-1,1)$, $\rho_0\ge 0$ belongs to $L^3((0,1))$ while $\varphi^\pm\ge 0$ is in 
$W^{1,\infty}([0,T])$;
\item for $\a=1$, same assumptions except that $\varphi^\pm>0$.
\end{enumerate}
\end{Def}

The notion of admissible weak solution of the nonlinear diffusion equations (\ref{NLDiff}) 
or (\ref{NLDiff-1}) being less obvious than  in the case of the Carleman model, we recall
it below.

\begin{Def}\lb{DEF-ADMIS-SOL}
An admissible weak solution of $(\ref{NLDiff})$ or $(\ref{NLDiff-1})$ with Dirichlet boundary 
conditions
\be
\rho(t,0)=2\varphi^-\,,\quad\rho(t,1)=2\varphi^+\,,\qquad t\in(0,T)
\ee
and initial condition
\be
\rho(0,x)=\rho_0(x)\,,\qquad x\in(0,1)
\ee
for admissible data $(\rho_0,\varphi^\pm)$ is a nonnegative element of $L^2((0,T)\times(0,1))$
such that $\d_x\rho^{1-\a}\in L^2((0,T)\times(0,1))$ if $\a\not=1$ while 
$\ln\rho\in L^1_{loc}( (0,T)\times(0,1))$ with $\d_x\ln\rho\in L^2((0,T)\times(0,1))$ if $\a=1$,
that satisfies $(\ref{BCDiff})$ for a.e. $t\in(0,T)$, together with the identities
\be
\int_0^T\int_0^1\left(\rho\d_t\phi-\tfrac{1}{2(1-\a)}\d_x\rho^{1-\a}\d_x\phi\right)(t,x)dxdt
+\int_0^1\rho_0(x)\phi(x,0)dxdt=0\,,
\end{equation}
for all $\phi\in C([0,T]\times[0,1]) \cap H^1([0,T]\times[0,1])$ vanishing in $x=0$, $x=1$ and 
$t=T$ if $\a\not=1$, whereas, for $\a=1$,
\be
\int_0^T\int_0^1\left(\rho\d_t\phi-\tfrac{1}{2}\d_x\ln\rho\d_x\phi\right)(t,x)dxdt
+\int_0^1\rho_0(x)\phi(x,0)dxdt=0\,.
\end{equation}
\end{Def}

Notice that, in what follows, we only need to know that the boundary value problem for (\ref{NLDiff})  or
(\ref{NLDiff-1}) has at most one admissible solution, as the existence of such solutions 
will result from that of a solution of the generalized Carleman model in the hydrodynamic 
limit. The proof of the following lemma is classical:

\begin{Lem}\lb{LM-UNIQ-DIFF}
For each set of admissible data $(\rho_0,\varphi^\pm)$, the nonlinear diffusion equation
$(\ref{NLDiff})$ for $-1\le\a<1$ or $(\ref{NLDiff-1})$ if $\a=1$ has at most one admissible 
solution.
\end{Lem}

\subsection{The convergence result.}

The main result in the present paper is the following convergence theorem.

\begin{Thm}
Let $(u^{in},v^{in},\varphi^\pm)$ be admissible data for the generalized Carleman system
on the time interval $[0,T]$. For each $\eps>0$, let $(u_\eps,v_\eps)$ be the solution
of the scaled Carleman system $(\ref{ScldCrlm})$ with initial condition $(\ref{InCond})$ and
boundary data $(\ref{BoundCond})$. 

Then, in the limit as $\eps\to 0$, the macroscopic density 
$$
\rho_\eps=u_\eps+v_\eps\to\rho\quad\hbox{ in }L^2_{loc}([0,T]\times[0,1])
$$
where $\rho$ is the generalized solution of $(\ref{NLDiff})$ if $\a\in[-1,1)$ or of $(\ref{NLDiff-1})$
if $\a=1$, with initial and boundary conditions
$$
\rho(x,0)= u^{in}(x)+v^{in}(x)\,,\qquad x\in(0,1)\,,
$$
$$
\rho(t,0)=2\varphi^-\,,\quad\rho(t,1)=2\varphi^+\,,\qquad t\in(0,T)\,.
$$
\end{Thm}

\section{Uniform bounds on the density and current}

As mentioned in the introduction, the uniform (in $\eps$) bounds
on the density and current are consequences of the equation that
is satisfied by some notion of relative entropy, which we shall
introduce below.

From now on, we set
\be\lb{Def-nu}
\nu=\max(\|\vphi^+\|_{W^{1,\infty}(0,T)},
    \|\vphi^-\|_{W^{1,\infty}(0,T)})
\ee
and 
\be\lb{Def-phim}
\varphi_m=\inf\{\varphi^\pm(t)\,|\,0\le t\le T\}>0\,.
\ee

We begin with some elementary background on convex functions.

\subsection{A family of convex functions}

Let $\phi:\,\bR_+\to\bR$ be a convex function that is $C^2$ on
$\bR_+^*$, and that satisfies
$$
\phi'_d(0)\le 0
$$
and
$$
\frac{\phi(y)}{y}\to+\infty\hbox{ as }y\to+\infty\,.
$$
Define the Legendre dual $\phi^*$ of $\phi$ by
$$
\phi^*(\xi)=\sup_{y\ge 0}(\xi y-\phi(y))\,.
$$
This definition clearly implies that
$$
\xi y\le\phi(y)+\phi^*(\xi)\hbox{ for each }y,\xi\in\bR_+\,.
$$
This inequality is a generalization of Young's classical
inequality for
$$
\phi(y)=\frac{y^p}{p}\hbox{ and }\phi^*(\xi)=\frac{\xi^{p'}}{p'}
$$
where $p$ and $p'$ are dual exponents in the sense of H\"older's
inequality:
$$
p\in[1,\infty)\hbox{ and }p'=\frac{p}{p-1}\,.
$$

More specifically, we shall use the family of convex functions
below:

\smallskip
a) for $\b\in[-1,1)$, the function $\phi_\b$ is given by
$$
\phi_\b(y)=\tfrac1{2-\b}y^{2-\b}\,,\qquad y\ge 0\,;
$$

b) for $\b=1$, the function $\phi_1$ is given by
$$
\phi_1(y)=y\log y\hbox{ and }\phi_1(0)=0\,,\qquad y>0\,.
$$

\smallskip
Straightforward computations show that
$$
\phi_\b^*(\xi)=\tfrac{1-\b}{2-\b}\xi^{\frac{2-\b}{1-\b}}
$$
while
$$
\phi_1^*(\xi)=e^{\xi-1}
$$
for each $\xi>0$ and $\b\in[-1,1)$.

With $\nu$ as in (\ref{Def-nu}), for each $\b\in[-1,1]$, we set
\be\lb{Def-Phib}
\Phi_\b(y)=\phi_\b(y)+\phi_\b^*(\phi'_\b(\nu)+1)\,,\quad y\ge 0\,.
\ee
Young's inequality entails
$$
\Phi_\b(y)\ge(\phi'_\b(\nu)+1)y\,,\quad y\ge0\,,
$$
and since $\Phi'_\b=\phi'_\b$, one has
\be\lb{CsqcYoung}
\Phi_\b(y)-\Phi'_\b(\nu)y\ge y\,,\quad \b\in[-1,1]\,,\,\,y\ge 0\,.
\ee

\subsection{The relative entropy and entropy production rate}

Our bounds on the Carleman system involve the notion of relative
entropy with respect to some suitable profile that satisfies the
same boundary conditions as the solution of the limiting diffusion
equation. There are many possible ways of choosing this profile.
When the boundary data are constant, a logical choice would be to
pick the stationary solution of the diffusion equation with those
boundary data. In the case of time dependent boundary conditions,
this choice is less natural, and we simply define this profile
to be the convex combination of boundary data:
$$
f(t,x):=(1-x)\vphi^-(t)+x\vphi^+(t)\,.
$$

Given $\Phi:\,\bR_+\to\bR$, a $C^2$ convex function, we define the
relative entropy of the 2-velocity density $(u,v)\equiv (u(x),v(x))
\in(\bR_+)^2$ to be
\be\lb{DefRelEntr}
H_\b[u,v|f]=
\int_0^1[\Phi_\b(u)+\Phi_\b(v)-2\Phi_\b(f)-\Phi'_\b(f)(u+v-2f)](x)dx\,,
\ee
and the entropy production rate as
\be\lb{DefEntrPrd}
P_\b[u,v]=\int_0^1(\Phi'_\b(u)-\Phi'_\b(v))(u+v)^\a(u-v)(x)dx\ge 0\,.
\ee

Assuming that $(u_\eps,v_\eps)$ is a solution of the scaled Carleman
problem (\ref{ScldCrlm}), we next determine the equation governing
the evolution of $H_\b[u_\eps,v_\eps|f]$. Multiplying each side of
the first equation in (\ref{ScldCrlm}) by $\Phi'_\b(u_\eps)$ and each
side of the second equation in (\ref{ScldCrlm}) by $\Phi'_\b(v_\eps)$,
one finds, upon adding both resulting equalities and integrating in
$x$ over $[0,1]$:
$$
\ba
\frac{d}{dt}\int_0^1\left(\Phi_\b(u_\eps)+\Phi_\b(v_\eps)\right)(t,x)dx
+
\frac1\eps\int_0^1\d_x\left(\Phi_\b(u_\eps)-\Phi_\b(v_\eps)\right)(t,x)dx
\\
=
-\frac1{\eps^2}P_\b[u_\eps,v_\eps]\,.
\ea
$$
It is actually more convenient to rearrange this equality as follows:
\be\lb{EntrId}
\ba
\frac{d}{dt}\int_0^1\left(\Phi_\b(u_\eps)+\Phi_\b(v_\eps)
    -\Phi'_\b(f)(u_\eps+v_\eps)\right)(t,x)dx
\\
+\frac1\eps\int_0^1\d_x\left(\Phi_\b(u_\eps)-\Phi_\b(v_\eps)
    -\Phi'_\b(f)(u_\eps-v_\eps)\right)(t,x)dx
\\
=-\frac1{\eps^2}P_\b[u_\eps,v_\eps]
-\frac1\eps\int_0^1\d_x\Phi_\b(f)(u_\eps-v_\eps)(t,x)dx
\\
+\int_0^1\d_t\Phi_\b(f)(u_\eps+v_\eps)(t,x)dx\,.
\ea
\ee
Observe that
\be\lb{IntBdrTrm}
\ba
\frac1\eps\int_0^1\d_x\left(\Phi_\b(u_\eps)-\Phi_\b(v_\eps)
   -\Phi_\b(f)(u_\eps-v_\eps)\right)(t,x)dx
\\
=
\frac1\eps[\Phi_\b(u_\eps)-\Phi_\b(f)-\Phi'_\b(f)(u_\eps-f)]_0^1
\\
-
\frac1\eps[\Phi_\b(v_\eps)-\Phi_\b(f)-\Phi'_\b(f)(v_\eps-f)]_0^1
\\
=
\frac1\eps\left(\Phi_\b(u_\eps)-\Phi_\b(f)-\Phi'_\b(f)(u_\eps-f)\right)(t,1)
\\
+
\frac1\eps\left(\Phi_\b(v_\eps)-\Phi_\b(f)-\Phi'_\b(f)(v_\eps-f)\right)(t,0)
\ea
\ee
since
$$
u_\eps(t,0)=f(t,0)\hbox{ and }v_\eps(t,1)=f(t,1)\,.
$$
By convexity of $\Phi_\b$, both terms in the last right hand side
of (\ref{IntBdrTrm}) are nonnegative; hence
\be\lb{SgnIntBdrTrm}
\frac1\eps\int_0^1\d_x\left(\Phi_\b(u_\eps)-\Phi_\b(v_\eps)
   -\Phi'_\b(f)(u_\eps-v_\eps)\right)(t,x)dx\ge 0\,.
\ee

Next, we formulate the equality (\ref{EntrId}) in terms of the
relative entropy $H_\b[u_\eps,v_\eps|f]$ as follows:
\be\lb{RltEntr<}
\ba
\frac{d}{dt}H_\b[u_\eps,v_\eps|f]+\frac1{\eps^2}P_\b[u_\eps,v_\eps]
\le
-\frac1\eps\int_0^1\d_x\Phi_\b(f)(u_\eps-v_\eps)(t,x)dx
\\
+\int_0^1\d_t\Phi'_\b(f)(u_\eps+v_\eps)(t,x)dx
+\int_0^1\Phi''_\b(f)\d_tf^2(t,x)dx\,.
\ea
\ee
We estimate the first term on the right hand side of the
inequality (\ref{RltEntr<}) in terms of the current
$j_\eps=\frac1\eps(u_\eps-v_\eps)$, as follows: for each
$\g\in(0,1)$, one has
\be\lb{1stTrm}
-\frac1\eps\int_0^1\d_x\Phi_\b(f)(u_\eps-v_\eps)(t,x)dx
\le
\tfrac{\g}2\int_0^1j_\eps(t,x)^2dx
+
\tfrac1{2\g}\|\d_x\Phi_\b(f)(t,x)\|^2_{L^\infty}\,.
\ee
Then, the inequality (\ref{CsqcYoung}) satisfied by $\Phi_\b$
implies that
\be\lb{CsqcYoung2}
\Phi_\b(y)-\Phi'_\b(f(t,x))y\ge\Phi_\b(y)-\Phi'_\b(\nu)y\ge y
\ee
for each $y\ge 0$. Indeed
$$
f(t,x)\le\max(\|\vphi^-\|_{L^\infty(0,T)},\|\vphi^+\|_{L^\infty(0,T)})
$$
so that
$$
\Phi'_\b(f(t,x))\le\Phi'_\b(\nu)
$$
since $\Phi'_\b$ is nondecreasing, $\Phi_\b$ being convex.
Because of (\ref{CsqcYoung2}), one has
\be\lb{2ndTrm}
\ba
\int_0^1\d_t\Phi_\b(f)&(u_\eps+v_\eps)(t,x)dx
\\
&\le
\int_0^1|\d_t\Phi'_\b(f)|[\Phi_\b(u_\eps)+\Phi_\b(v_\eps)
    -\Phi'_\b(f)(u_\eps+v_\eps)](t,x)dx
\\
&\le
\|\d_t\Phi'_\b(f)\|_{L^\infty}(H_\b[u_\eps,v_\eps|f]
    +2\|\Phi_\b(f)-\Phi'_\b(f)f\|_{L^\infty})\,.
\ea
\ee
With this estimate for the second term on the right hand side of
(\ref{RltEntr<}), we recast this inequality as
\be\lb{RltEntr<2}
\ba
\frac{d}{dt}H_\b[u_\eps,v_\eps|f]+\frac1{\eps^2}P_\b[u_\eps,v_\eps]
\le
\tfrac{\g}2\int_0^1j_\eps(t,x)^2dx+\|\d_t\Phi'_\b(f)\|_{L^\infty}H_\b[u_\eps,v_\eps|f]
\\
+\tfrac1{2\g}\|\d_x\Phi_\b(f)\|^2_{L^\infty}+\|\Phi''_\b(f)\d_tf^2\|_{L^\infty}
\\
+2\|\d_t\Phi'_\b(f)\|_{L^\infty}\|\Phi_\b(f)-\Phi'_\b(f)f\|_{L^\infty}\,.
\ea
\ee

\subsection{Current estimate}

Our first objective is an $L^2_{t,x}$ bound on the current
$j_\eps=\frac1\eps(u_\eps-v_\eps)$.

Start from the relative entropy inequality (\ref{RltEntr<2}) with
the particular choice $\b=\a$ (we recall that $\a$ is the exponent
in the nonlinearity of the Carleman system (\ref{ScldCrlm})).

We are going to show that, for $\g$ small enough, the term
$$
\tfrac{\g}2\int_0^1j_\eps(t,x)^2dx
\hbox{ is dominated by }
\frac1{\eps^2}P_\a[u_\eps,v_\eps]\,.
$$
Indeed, if $\a=1$
$$
\ba
\frac1{\eps^2}P_\a[u_\eps,v_\eps]
&=\frac1{\eps^2}\int_0^1(u_\eps^2-v_\eps^2)(\log u_\eps-\log v_\eps)dx
\\
&=
\int_0^1j_\eps^2(u_\eps+v_\eps)\frac{\log u_\eps-\log v_\eps}{u_\eps-v_\eps}dx
\\
&\ge
\int_0^1j_\eps^2\frac{u_\eps+v_\eps}{\max(u_\eps,v_\eps)}dx
\ge
\int_0^1j_\eps^2dx\,.
\ea
$$
in view of the elementary inequality
$$
\frac{\log a-\log b}{a-b}\ge\frac1{\max(a,b)}\,,\quad\hbox{ for each }a,b>0\,.
$$
For $\a\in[0,1)$
$$
\ba
\frac1{\eps^2}P_\a[u_\eps,v_\eps]
&=\frac1{\eps^2}\int_0^1(u_\eps-v_\eps)(u_\eps+v_\eps)^\a(u^{1-\a}_\eps-v^{1-\a}_\eps)dx
\\
&=
\int_0^1j_\eps^2(u_\eps+v_\eps)^\a\frac{u^{1-\a}_\eps-v^{1-\a}_\eps}{u_\eps-v_\eps}dx
\\
&\ge
(1-\a)\int_0^1j_\eps^2\frac{(u_\eps+v_\eps)^\a}{\max(u_\eps,v_\eps)^\a}dx
\ge
(1-\a)\int_0^1j_\eps^2dx\,.
\ea
$$
Finally, for $\a\in[-1,0)$, we separate $(0,1)$ into the region where
$u_\eps(t,x)>v_\eps(t,x)$ and its complement, where $v_\eps(t,x)\ge u_\eps(t,x)$.
$$
\ba
\frac1{\eps^2}P_\a[u_\eps,v_\eps]
&=\frac1{\eps^2}\int_0^1(u_\eps-v_\eps)(u_\eps+v_\eps)^\a(u^{1-\a}_\eps-v^{1-\a}_\eps)dx
\\
&=
\int_0^1j_\eps^2(u_\eps+v_\eps)^\a\frac{u^{1-\a}_\eps-v^{1-\a}_\eps}{u_\eps-v_\eps}dx
\\
&=\int_{u_\eps>v_\eps}
+\int_{v_\eps\ge u_\eps}j_\eps^2(u_\eps+v_\eps)^\a
\frac{u^{1-\a}_\eps-v^{1-\a}_\eps}{u_\eps-v_\eps}dx\,.
\ea
$$
The first integral is estimated as follows:
$$
\ba
\int_{u_\eps>u_\eps}j_\eps^2(u_\eps+v_\eps)^\a
\frac{u^{1-\a}_\eps-v^{1-\a}_\eps}{u_\eps-v_\eps}dx
&\ge
\int_{u_\eps>v_\eps}j_\eps^2(2u_\eps)^\a
\frac{u^{1-\a}_\eps-v^{1-\a}_\eps}{u_\eps-v_\eps}dx
\\
&=
2^\a\int_{u_\eps>v_\eps}j_\eps^2
\frac{u^{1+|\a|}_\eps-v^{1+|\a|}_\eps}{u_\eps^{1+|\a|}-u_\eps^{|\a|}v_\eps}dx
\\
&\ge
2^\a\int_{u_\eps>v_\eps}j_\eps^2dx
\ea
$$
The second integral is estimated similarly. In all cases,
\be\lb{Curr<EntrPrd}
\frac1{\eps^2}P_\a[u_\eps,v_\eps]\ge C_\a\int_{u_\eps>v_\eps}j_\eps^2dx
\ee
with
$$
C_1=1\hbox{ while }C_\a=1-\a\hbox{ if }0\le\a<1
\hbox{ and }C_\a=2^\a\hbox{ if }-1\le\a<0\,.
$$

\begin{Prop}\lb{PR-BNDJ}
Let $(u_\eps,v_\eps)$ be a solution of the scaled, generalized
Carleman system $(\ref{ScldCrlm})$ with admissible initial and
boundary data. Then there exists a positive constant $J\equiv
J(\nu,\vphi_m,\a,T,u^{in},v^{in})>0$ such that the current
$j_\eps=\frac1\eps(u_\eps-v_\eps)$ satisfies
$$
\int_0^T\int_0^1j^2_\eps(t,x)dxdt\le J
$$
for each $\eps>0$.
\end{Prop}

\begin{proof}
Start from (\ref{RltEntr<2}), and pick $\g=C_\a/2$. Upon integrating
both sides of (\ref{RltEntr<2}) over the time interval $[0,t]$, one
obtains
$$
\ba
H_\a[u_\eps,v_\eps|f](t)&+\tfrac{3C_\a}4\int_0^t\int_0^1j_\eps(s,x)^2dxds
\\
&\le\|\d_t\Phi'_\a(f)\|_{L^\infty}\int_0^tH_\a[u_\eps,v_\eps|f](s)ds
\\
&\qquad+M_0(\nu,\a,\vphi_m)t+H_\a[u^{in},v^{in}|f]\,,
\ea
$$
where
$$
\ba
M_0(\nu,\a,\vphi_m)\ge\tfrac1{C_\a}\|\d_x\Phi_\a(f)\|^2_{L^\infty}
+\|\Phi''_\a(f)\d_tf^2\|_{L^\infty}
\\
+2\|\d_t\Phi'_\a(f)\|_{L^\infty}\|\Phi_\a(f)-\Phi'_\a(f)f\|_{L^\infty}\,.
\ea
$$
Gronwall's inequality implies that
$$
\ba
\|\d_t\Phi'_\a(f)\|_{L^\infty}\int_0^tH_\a[u_\eps,v_\eps|f](s)ds
+M_0(\nu,\a,\vphi_m)t+H_\a[u^{in},v^{in}|f]
\\
\le
(M_0(\nu,\a,\vphi_m)t+H_\a[u^{in},v^{in}|f])e^{t\|\d_t\Phi'_\a(f)\|_{L^\infty}}
\ea
$$
from which we deduce that
$$
\int_0^t\int_0^1j_\eps(s,x)^2dxds\le\tfrac4{3C_\a}(M_0(\nu,\a,\vphi_m)t
+H_\a[u^{in},v^{in}|f])e^{t\|\d_t\Phi'_\a(f)\|_{L^\infty}}
$$
\end{proof}

\subsection{Density estimate}

With the current estimate at our disposal, we next obtain
an $L^2_{loc}(\bR_+;L^2_x)$ bound on the macroscopic
density $\rho_\eps=u_\eps+v_\eps$.

To do this, we apply (\ref{RltEntr<2}) with $\b=0$; here,
the entropy production rate is useless so that we actually
discard it from the left-hand side:
$$
\frac{d}{dt}H_0[u_\eps,v_\eps|f]\le\tfrac{\g}2\int_0^1j_\eps(t,x)^2dx
+\|\d_t\Phi'_0(f)\|_{L^\infty}H_0[u_\eps,v_\eps|f]+M_1(\nu,\g)
$$
where
$$
\ba
M_1(\nu,\g)&\ge
\tfrac1{2\g}\|\d_xf^2\|^2_{L^\infty}
+\|\d_tf^2\|_{L^\infty}+(1+\nu)^2\|\d_tf\|_{L^\infty}
\\
&\ge\tfrac1{2\g}\|\d_x\Phi_0(f)\|^2_{L^\infty}+\|\Phi''_0(f)\d_tf^2\|_{L^\infty}
\\
&+2\|\d_t\Phi'_0(f)\|_{L^\infty}\|\Phi_0(f)-\Phi'_0(f)f\|_{L^\infty}\,.
\ea
$$
Integrating the relative entropy inequality above over the
time interval $[0,t]$, we obtain
$$
\ba
H_0[u_\eps,v_\eps|f](t)\le H_0[u^{in},v^{in}|f]
    +\tfrac{\g}2J(\nu,\varphi_m,\a,t,u^{in},v^{in})+M_1(\nu,\g)t
\\
+\|\d_t\Phi'_0(f)\|_{L^\infty}\int_0^tH_0[u_\eps,v_\eps|f](s)ds\,.
\ea
$$
Gronwall's inequality implies that
$$
\ba
{}&H_0[u_\eps,v_\eps|f](t)
\\
&\le\left(H_0[u^{in},v^{in}|f]
+\tfrac{\g}2J(\nu,\vphi_m,\a,t,u^{in},v^{in})+M_1(\nu,\g)t\right)
e^{t\|\d_t\Phi'_0(f)\|_{L^\infty}}
\ea
$$
On the other hand,
$$
H_0[u_\eps,v_\eps|f](t)=
\tfrac12\left(\|u_\eps(t)-f(t)\|^2_{L^2_x}
+\|v_\eps(t)-f(t)\|^2_{L^2_x}\right)\,.
$$
Summarizing the discussion above, we have proved the following
density bound.

\begin{Prop}\lb{PR-BNDRHO}
Let $(u_\eps,v_\eps)$ be a solution of the scaled, generalized
Carleman system $(\ref{ScldCrlm})$ with admissible initial and
boundary data. Then there exists a positive constant $K\equiv
K(\nu,\vphi_m,\a,T,u^{in},v^{in})>0$ such that the macroscopic
density $\rho_\eps=u_\eps+v_\eps$ satisfies
$$
\int_0^T\int_0^1\rho^2_\eps(t,x)dxdt\le K
$$
for each $\eps>0$.
\end{Prop}

\subsection{An entropy production estimate}

We shall conclude this section with a further nonlinear estimate
that eventually controls some flux in the diffusion limit.

We distinguish the cases $\a\in[0,1]$ and $\a\in[-1,0)$.

If $\a\in[0,1]$, we use (\ref{RltEntr<2}) with $\b=\pm1$. 

Indeed, the entropy production rate in the case $\b=-1$ is
$$
\ba
P_{-1}[u_\eps,v_\eps]&=\int_0^1(u_\eps+v_\eps)^\a
(u_\eps-v_\eps)(\Phi'_{-1}(u_\eps)-\Phi'_{-1}(v_\eps))dx
\\
&=\int_0^1(u_\eps+v_\eps)^\a(u_\eps-v_\eps)(u^2_\eps-v^2_\eps)dx
=\int_0^1\rho_\eps^{1+\a}(\eps j_\eps)^2dx\,.
\ea
$$
Hence the relative entropy inequality (\ref{RltEntr<2}) with
$\b=-1$ becomes
$$
\ba
\frac{d}{dt}H_{-1}[u_\eps,v_\eps|f]
    +\int_0^1\rho_\eps^{1+\a}j_\eps^2dx
\le
\tfrac{\g}2\int_0^1j_\eps(t,x)^2dx+M_2(\nu,\g)
\\
+\|\d_t\Phi'_{-1}(f)\|_{L^\infty}H_{-1}[u_\eps,v_\eps|f]
\ea
$$
with
$$
\ba
M_2(\nu,\g)&\ge2\|\d_t\Phi'_{-1}(f)\|_{L^\infty}
    \|\Phi_{-1}(f)-\Phi'_{-1}(f)f\|_{L^\infty}
\\
&\quad+\tfrac1{2\g}\|\d_x\Phi_{-1}(f)\|^2_{L^\infty}
    +\|\Phi''_{-1}(f)\d_tf^2\|_{L^\infty}\,.
\ea
$$
Integrating the relative entropy inequality above in the
time variable, and applying Proposition \ref{PR-BNDJ}, we
obtain
$$
\ba
H_{-1}[u_\eps,v_\eps|f](t)
    +\int_0^t\int_0^1\rho_\eps^{1+\a}j_\eps^2dx
\le
\|\d_t\Phi'_{-1}(f)\|_{L^\infty}
    \int_0^tH_{-1}[u_\eps,v_\eps|f](s)ds
\\
+H_{-1}[u^{in},v^{in}|f]
+\tfrac{\g}2J(\nu,\vphi_m,\a,T,u^{in},v^{in})+M_2(\nu,\g)t
\ea
$$
so that, by Gronwall's inequality, we arrive at an estimate
of the form
$$
\ba
{}&\int_0^t\int_0^1\rho_\eps^{1+\a}j_\eps^2dx
\\
&\le
\left(H_{-1}[u^{in},v^{in}|f]
+\tfrac{\g}2J(\nu,\vphi_m,\a,T,u^{in},v^{in})+M_2(\nu,\g)t\right)
e^{t\|\d_t\Phi'_{-1}(f)\|_{L^\infty}}\,.
\ea
$$

Doing the same with $\b=1$ leads to
$$
\ba
P_{1}[u_\eps,v_\eps]&=\int_0^1(u_\eps+v_\eps)^\a
(u_\eps-v_\eps)(\ln(u_\eps)-\ln(v_\eps))dx
\\
&\ge\int_0^1(u_\eps+v_\eps)^\a(u_\eps-v_\eps)\frac{u_\eps-v_\eps}{u_\eps-v_\eps}dx
=\int_0^1\rho_\eps^{\a-1}(\eps j_\eps)^2dx\,.
\ea
$$
The estimate above rests on the mean value theorem:
$$
(a-b)(\ln a-\ln b)\ge\frac{|a-b|}{\max(a,b)}\ge\frac{|a-b|}{a+b}\,,\qquad a,b>0\,.
$$
Integrating in $t$ the relative entropy inequality (\ref{RltEntr<2}) with
$\b=1$ and applying the Gronwall inequality as above gives
$$
\ba
{}&\int_0^t\int_0^1\rho_\eps^{\a-1}j_\eps^2dx
\\
&\le
\left(H_{-1}[u^{in},v^{in}|f]
+\tfrac{\g}2J(\nu,\vphi_m,\a,T,u^{in},v^{in})+M_3(\nu,\vphi_m,\g)t\right)
e^{t\|\d_t\Phi'_{-1}(f)\|_{L^\infty}}\,,
\ea
$$
where
$$
\ba
M_3(\nu,\vphi_m,\g)&\ge2\|\d_t\Phi'_{1}(f)\|_{L^\infty}
    \|\Phi_{1}(f)-\Phi'_{1}(f)f\|_{L^\infty}
\\
&\quad+\tfrac1{2\g}\|\d_x\Phi_{1}(f)\|^2_{L^\infty}
    +\|\Phi''_{1}(f)\d_tf^2\|_{L^\infty}\,.
\ea
$$

If on the other hand $\a\in[-1,0)$, we use (\ref{RltEntr<2})
with $\b=|\a|$. Notice that $\Phi''_{|\a|}$ is decreasing on
$\bR_+^*$: indeed
$$
\ba
\Phi''_{|\a|}(z)&=(1-|\a|)z^{-|\a|}&&\hbox{ if }\a\in(-1,0)\,,
\\
\Phi''_{|\a|}(z)&=z^{-1}&&\hbox{ if }\a=-1\,.
\ea
$$
Hence, for each $a,b>0$, the mean value theorem implies that,
for some $\th\in(0,1)$
$$
(a-b)(\Phi'_{|\a|}(a)-\Phi'_{|\a|}(b))
=(a-b)^2\Phi''_{|\a|}((1-\th)a+\th b)
\ge(a-b)^2\Phi''_{|\a|}(a+b)\,.
$$
Therefore
$$
\ba
P_{|\a|}[u_\eps,v_\eps]&=\int_0^1(u_\eps+v_\eps)^\a
(u_\eps-v_\eps)(\Phi'_{|\a|}(u_\eps)-\Phi'_{|\a|}(v_\eps))dx
\\
&\ge\int_0^1(u_\eps-v_\eps)^2(u_\eps+v_\eps)^\a \Phi''_{|\a|}(u_\eps+v_\eps)dx
=C'_\a\int_0^1(\eps j_\eps)^2\rho_\eps^{2\a}dx\,,
\ea
$$
with
$$
C'_\a=1\hbox{ if }\a=-1\hbox{ and }C_\a=1+\a\hbox{ if }\a\in(-1,0)\,.
$$
The relative entropy inequality (\ref{RltEntr<2}) with $\b=|\a|$
becomes
$$
\ba
\frac{d}{dt}H_{|\a|}[u_\eps,v_\eps|f]
    +C_\a\int_0^1\rho_\eps^{2\a}j_\eps^2dx
\le
\tfrac{\g}2\int_0^1j_\eps(t,x)^2dx+M_4(\nu,\vphi_m,\g)
\\
+\|\d_t\Phi'_{|\a|}(f)\|_{L^\infty}H_{|\a|}[u_\eps,v_\eps|f]
\ea
$$
with
$$
\ba
M_4(\nu,\vphi_m,\g)&\ge2\|\d_t\Phi'_{|\a|}(f)\|_{L^\infty}
    \|\Phi_{|\a|}(f)-\Phi'_{|\a|}(f)f\|_{L^\infty}
\\
&\quad+\tfrac1{2\g}\|\d_x\Phi_{|\a|}(f)\|^2_{L^\infty}
    +\|\Phi''_{|\a|}(f)\d_tf^2\|_{L^\infty}\,.
\ea
$$
Integrating this inequality over the time interval $[0,t]$ and
applying Gronwall's inequality and Proposition \ref{PR-BNDJ} as
before, we obtain
$$
\ba
{}&C_\a\int_0^t\int_0^1\rho_\eps^{2\a}j_\eps^2dx
\\
&\le
\left(H_{|\a|}[u^{in},v^{in}|f]
+\tfrac{\g}2J(\nu,\varphi_m,\a,T,u^{in},v^{in})+M_4(\nu,\vphi_m,\g)t\right)
e^{t\|\d_t\Phi'_{|\a|}(f)\|_{L^\infty}}\,.
\ea
$$

Summarizing, we have proved

\begin{Prop}\lb{PR-BNDRHOJ}
Let $(u_\eps,v_\eps)$ be a solution of the scaled, generalized
Carleman system $(\ref{ScldCrlm})$ with admissible initial and
boundary data. Then there exists a positive constant $L\equiv
L(\nu,\vphi_m,\a,t,u^{in},v^{in})>0$ such that
$$
\ba
\int_0^T\int_0^1\rho^{\a+1}_\eps j^2_\eps(t,x)dxdt\le L\,,
\\
\int_0^T\int_0^1\rho^{\a-1}_\eps j^2_\eps(t,x)dxdt\le L\,,
\ea
$$
for each $T>0$ and $\eps>0$ if $\a\in[0,1]$, and
$$
\int_0^T\int_0^1\rho^{2\a}_\eps j^2_\eps(t,x)dxdt\le L
$$
whenever $\a\in[-1,1]$.
\end{Prop}

In the case $\a\in[0,1]$, the second bound follows from the
first and Proposition \ref{PR-BNDJ} by H\"older's inequality.

\section{The nonlinear diffusion limit}

Let $(u^{in},v^{in})$, $\varphi^+$ and $\varphi^-$ be admissible
initial and boundary data; then, for each $\eps>0$, let
$(u_\eps,v_\eps)$ be the solution to the scaled, generalized
Carleman system (\ref{ScldCrlm}).

It follows from Propositions \ref{PR-BNDJ} and \ref{PR-BNDRHO} 
that, for each $T>0$, one has
$$\
\ba
\|\rho_\eps\|_{L^\infty(0,T;L^2(0,1))}
    &\le K(\nu,\vphi_m,\a,T,u^{in},v^{in})
\hbox{ and }
\\
\|j_\eps\|_{L^2([0,T]\times[0,1])}
    &\le J(\nu,\vphi_m,\a,T,u^{in},v^{in})
\ea
$$
for each $\eps>0$. By the Banach-Alaoglu theorem, for each
$T>0$
\be\lb{WCp-rho}
\hbox{the family }\rho_\eps\hbox{ is relatively compact in }
    L^\infty(0,T;L^2(0,1))\hbox{ weak-*}
\ee
while
\be\lb{WCp-j}
\hbox{the family }j_\eps\hbox{ is relatively compact in }
    L^2([0,T]\times[0,1])\hbox{ weak.}
\ee

Summing both equations in Carleman's system implies that
$$
\d_t\rho_\eps=-\d_xj_\eps
$$
so that
$$
\d_t\rho_\eps\hbox{ is bounded in }L^2(0,T;H^{-1}(0,1))\,;
$$
with the bound on $\rho_\eps$ coming from Proposition
\ref{PR-BNDRHO}, this control implies that
\be\lb{SCp-rho-t}
\rho_\eps\hbox{ is relatively compact in }
C([0,T];H^{-1}(0,1))
\ee
by Arzela-Ascoli's theorem.

Likewise, since
$$
u_\eps=\frac{\rho_\eps+\eps j_\eps}2
	\hbox{ and }v_\eps=\frac{\rho_\eps-\eps j_\eps}2 \,,
$$
it follows that
$$
\ba
\d_x(\rho_\eps+\eps j_\eps)=-2\rho_\eps^\a j_\eps-\eps\d_t(\rho_\eps+\eps j_\eps)
\hbox{Ê is bounded in }L^2([0,1];H^{-1}(0,T))
\\
\d_x(\rho_\eps-\eps j_\eps)=-2\rho_\eps^\a j_\eps+\eps\d_t(\rho_\eps-\eps j_\eps)
\hbox{Ê is bounded in }L^2([0,1];H^{-1}(0,T))
\ea
$$
because of Propositions \ref{PR-BNDJ}, \ref{PR-BNDRHO} and
\ref{PR-BNDRHOJ}.
Hence
\be\lb{SCp-rho-x}
\rho_\eps\pm\eps j_\eps\hbox{ is relatively compact in }
	C([0,1];H^{-1}(0,T))
\ee
by Arzela-Ascoli's theorem.

Because on the nonlinearities that are present both in
the scaled Carleman system (\ref{ScldCrlm}) and in the
limiting nonlinear diffusion equation (\ref{NLDiff}),
weak compactness results as above are not enough to pass
to the limit as $\eps\to 0$. Strong $L^2$ compactness of
the family $\rho_\eps$ is obtained by the same argument
as in \cite{MarcatiMilani} (see also \cite{LionsToscani}), which we recall below.

Consider the vector fields
$$
p_\eps=(\rho_\eps,j_\eps)\hbox{ and }
q_\eps=(-\rho_\eps,\eps^2j_\eps)\,.
$$
By Propositions \ref{PR-BNDJ} and \ref{PR-BNDRHO}, both
vector fields satisfy
$$
p_\eps\hbox{ and }q_\eps\hbox{ are bounded in }
L^2([0,T]\times[0,1];\bR^2)\,.
$$
Summing the first two equations in the scaled Carleman
system (\ref{ScldCrlm}) shows that
$$
\Div_{t,x}p_\eps=\d_t\rho_\eps+\d_xj_\eps=0
$$
while
$$
\Curl_{t,x}q_\eps=\eps^2\d_tj_\eps+\d_x\rho_\eps
    =-2\rho_\eps^\a j_\eps\,.
$$
By Proposition \ref{PR-BNDRHOJ}, we therefore have
$$
\Div_{t,x}p_\eps\hbox{ and }\Curl_{t,x}q_\eps
\hbox{ bounded in }L^2([0,T]\times[0,1])\,.
$$
Pick any sequence $\eps_n\to 0$ such that
$$
\rho_{\eps_n}\wto\rho\hbox{ and }
j_{\eps_n}\wto j\hbox{ in }L^2([0,T]\times[0,1])
$$
as $n\to\infty$. By compensated compactness (the div-curl
lemma in \cite{Murat}), we find that
$$
p_{\eps_n}\cdot q_{\eps_n}
=-\rho^2_{\eps_n}+\eps_n^2j^2_{\eps_n}
\wto
p\cdot q=\rho^2
$$
in the sense of Radon measures on $(0,T)\times(0,1)$
as $\eps_n\to 0$. Because of Proposition \ref{PR-BNDJ},
this implies that
$$
\rho_{\eps_n}^2\wto\rho^2
    \hbox{ in the sense of Radon measures}
$$
which implies in turn that the family
\be\lb{SCp-rho}
\rho_\eps\hbox{ is relatively compact in }
L^2([0,T]\times[0,1])\hbox{ strong.}
\ee

Let then $(\rho,j,q)$ be a weak limit point in $L^2([0,T]\times[0,1])$
of the family $(\rho_\eps,j_\eps,\rho_\eps^\a j_\eps)$ as $\eps\to 0$ 
--- the existence of such limits points being guaranteed by  the bounds 
in Propositions \ref{PR-BNDJ}, \ref{PR-BNDRHO} and \ref{PR-BNDRHOJ} 
together with the Banach-Alaoglu Theorem.

Because of (\ref{SCp-rho-t}) and of the initial condition of the
Carleman system (\ref{ScldCrlm}), one has
\be\lb{LimCondin}
\rho\in C([0,T];H^{-1}(0,1))\hbox{ and }\rho\rstr_{t=0}=u^{in}+v^{in}\,.
\ee

Because of (\ref{SCp-rho-x}) and of the boundary conditions of the
scaled Carleman system (\ref{ScldCrlm}), one has
\be\lb{LimBndCond}
\rho\in C([0,1];H^{-1}(0,T))\hbox{ and }\rho\rstr_{x=0}=2\varphi^-\,,
	\hbox{Êwhile }\rho\rstr_{x=1}=2\varphi^+\,.
\ee

Summing both equations in the scaled Carleman system shows that
$$
\d_t\rho_\eps+\d_xj_\eps=0
$$
so that, by passing to the limit in the relation above, one arrives at
\be\lb{LmContEq}
\d_t\rho+\d_xj=0
\ee
in the sense of distributions on $\bR_+^*\times(0,1)$.

Subtracting the second equation from the first in the scaled Carleman
system shows that
$$
\eps^2\d_tj_\eps+\d_x\rho_\eps=-2\rho^\a_\eps j_\eps
$$
so that, by passing to the limit in the equation above, one finds that
$$
\d_x\rho=-2q\,.
$$

Assume first that $\a\in[-1,0]$. By (\ref{SCp-rho}) and the weak-strong 
continuity of the product of two functions,\footnote{ For $f \ : \ \bR_+ 
\to \bR$ continuous and sublinear at infinity, if
$\rho_\eps \to \rho$ in $L^2$ and $j_\eps \rightharpoonup j$ in $L^2$,
$f(\rho_\eps)j_\eps \to f(\rho)j$ in $\mathcal{D'}$.}
$$
\rho^{|\a|}q=j
$$
so that
$$
j=-\tfrac12\rho^{|\a|}\d_x\rho=-\tfrac1{2(1-\a)}\d_x\rho^{1-\a}\,.
$$
Substituting this in (\ref{LmContEq}), we see that $\rho$ satisfies
\be\lb{LimDiffEq}
\ba
{}&\d_t\rho-\tfrac1{2(1-\a)}\d_x^2\rho^{1-\a}=0\,,
\\
&\rho(t,0)=2\varphi^-(t)\,,
\\
&\rho(t,1)=2\varphi^+(t)\,,
\\
&\rho(0,x)=u^{in}(x)+v^{in}(x)\,.
\ea
\ee
Since this problem has a unique solution, the whole sequence
$\rho_\eps\to\rho$ in $L^2([0,T]\times[0,1])$ as $\eps\to 0$, by the
compactness statement in (\ref{SCp-rho}).

The case $\a\in(0,1)$ is slightly more delicate. Since $\rho\ge 0$
a.e. but may vanish, we set $\eta>0$ and, using as above the
strong compactness statement (\ref{SCp-rho}) together with the
weak-strong continuity of the product of two functions,  
$$
\frac1{(\rho+\eta)^\a}q=\frac{\rho^\a}{(\rho+\eta)^\a}j
$$
so that
$$
\frac{\rho^\a}{(\rho+\eta)^\a}j=-\tfrac1{2(1-\a)}\d_x(\eta+\rho)^{1-\a}\,.
$$
We recall the argument in \cite{LionsToscani}. Letting $\eta\to 0$ 
shows that
$$
\indc_{\rho>0}j=-\tfrac1{2(1-\a)}\d_x\rho^{1-\a}\,.
$$
On the other hand, the second bound in Proposition \ref{PR-BNDRHOJ}
implies that
$$
\int_0^T\int_0^1\rho^{\a-1}j^2dxdt<+\infty
$$
so that
$$
j=0\hbox{ a.e. wherever }\rho=0\,.
$$
Hence
$$
\indc_{\rho>0}j=j
$$
so that the limiting nonlinear diffusion equation (\ref{LimDiffEq}) 
also holds in the case $\a\in(0,1)$.

It only remains to treat the case $\a=1$. Proceeding as above,
we obtain instead the relation
$$
\frac{\rho}{\rho+\eta}j=-\tfrac12\d_x\ln(\eta+\rho)\,.
$$
Integrating this relation in $x$ shows that
$$
\ln(\eta+\rho)(t,x)=\ln(\eta+\varphi^+(t))-2\int_0^x\frac{\rho}{\rho+\eta}j(t,z)dz
$$
so that, by Proposition \ref{PR-BNDJ} together with the bound 
$\varphi_m\le\varphi^+\le\nu$, we see that
$$
\int_0^T\int_0^1\sup_{\eta}|\ln(\eta+\rho)(t,x)|dtdx<+\infty\,.
$$
Hence, by dominated convergence, $\ln\rho\in L^1([0,T]\times[0,1])$
so that $\rho>0$ a.e. on $[0,T]\times[0,1]$. Therefore
$$
\frac{\rho}{\rho+\eta}j\to j\hbox{ in }L^2([0,T]\times[0,1])
$$
and
$$
\ln(\eta+\rho)\to\ln\rho\hbox{ in }L^1([0,T]\times[0,1])
$$
as $\eta\to 0$, so that
$$
j=-\tfrac12\d_x\ln\rho\,.
$$
Substituting this in (\ref{LmContEq}), we see that $\rho$ satisfies
\be\lb{LimDiffEq-ln}
\ba
{}&\d_t\rho-\tfrac12\d_x^2\ln\rho=0\,,
\\
&\rho(t,0)=2\varphi^-(t)\,,
\\
&\rho(t,1)=2\varphi^+(t)\,,
\\
&\rho(0,x)=u^{in}(x)+v^{in}(x)\,.
\ea
\ee
Again, since this problem has at most one solution, the whole sequence
$\rho_\eps$ converges to that solution $\rho$ in $L^2([0,T]\times[0,1])$ 
because of the compactness statement (\ref{SCp-rho}).

Thus we have established the nonlinear diffusion limit for each $\a\in[-1,1]$.



\begin{thebibliography}{99}

\bibitem{Carl}
{\sc T. Carleman.}
\newblock Probl\` emes  math\' ematiques dans la th\' eorie
cin\' etique des gaz,
\newblock Almqvist-Wiksells, Uppsala, 1957.

\bibitem{DoMa}
{\sc D. Donatelli, P. Marcati.}
\newblock Convergence of singular
limits for multi-D semilinear hyperbolic systems to parabolic systems,
\newblock {\em Trans. Amer. Math. Soc.} {\bf 356} (2004) 2093--2121.


\bibitem{Fitzg}
{\sc W. E. Fitzgibbon.}
\newblock The fluid-dynamical limit of the Carleman
equation with reflecting boundary,
\newblock {\em J. Nonlin. Anal. Th. Meth. Appl.} {\bf 6} (1982) 695--702.


\bibitem{Fitzg2}
{\sc W. E. Fitzgibbon.}
\newblock Initial-boundary value problem for the Carleman equation,
\newblock {\em Comput. Math. Appl.} {\bf 9} (1983) 519--525.


\bibitem{Gol}
{\sc S. Goldstein.}
\newblock On diffusion by discontinuous movements, and on the telegraph equation,
\newblock {\em Quart. J. Mech. Appl. Math.} {\bf 4} (1951) 129--156.

\bibitem{GolSal}
{\sc F. Golse, F. Salvarani.}
\newblock Work in preparation.



\bibitem{Kolod}
{\sc I. Kolodner.}
\newblock 
On Carleman's model for the Boltzmann equation and its generalizations,
\newblock {\em Ann. Mat. Pura Appl.} {\bf 63} (1963) 11--32.


\bibitem{Kurtz}
{\sc T. G. Kurtz.}
\newblock Convergence of sequences of semigroups of nonlinear
operators with an application to gas kinetics,
\newblock {\em Trans. Amer. Math. Soc.} {\bf 186} (1973) 259--272.


\bibitem{LionsToscani}
{\sc P. L. Lions, G. Toscani.}
\newblock Diffusive limits for finite velocities Boltzmann
kinetic models,
\newblock {\em Rev. Mat. Iberoamericana}, {\bf 13} (1997) 473--513.

\bibitem{MarcatiMilani}
{\sc P. Marcati, A. J. Milani.}
\newblock The one-dimensional Darcy's law as the limit of a compressible Euler
flow.
\newblock {\em J. Diff. Eq.} {\bf 84} (1990) 129--147.


\bibitem{MaRu}
{\sc P. Marcati, B. Rubino.}
\newblock Hyperbolic to parabolic relaxation theory for
quasilinear first order systems,
\newblock {\em J. Differential Equations}, {\bf 162}
(2000) 359--399.


\bibitem{McKean}
{\sc H. P. McKean.}
\newblock The central limit theorem for Carleman's equation,
\newblock {\em {Israel J. Math.}} {\bf 21} (1975) 54--92.

\bibitem{Mihalas}
{\sc D. Mihalas, B. Weibel-Mihalas.}
\newblock Foundations of Radiation Hydrodynamics
\newblock Dover, Mineola NY, 1999.


\bibitem{Murat}
{\sc F. Murat.}
\newblock Compacit\' e par compensation,
\newblock {\em {Ann. Scuola Norm. Sup. Pisa Cl. Sci.}} {\bf 5} (1978) 489--507.


\bibitem{PuTo}
{\sc A. Pulvirenti, G. Toscani.}
\newblock Fast diffusion as a limit of a two-velocity kinetic model,
\newblock {\em Rend. Circ. Mat. Palermo Suppl.} {\bf 45} part II (1996)
521--528.

\bibitem{ReedIllner1}
{\sc R. Illner, M.C. Reed.}
\newblock The decay of solutions of the Carleman model. 
\newblock Math. Methods Appl. Sci. {\bf 3} (1981), 121--127.

\bibitem{ReedIllner2}
{\sc R. Illner, M.C. Reed.}
\newblock Decay to equilibrium for the Carleman model in a box. 
\newblock SIAM J. Appl. Math. {\bf 44} (1984), 1067--1075.

\bibitem{Salv}
{\sc F. Salvarani.}
\newblock Diffusion limits for the initial-boundary value problem of the Goldstein-Taylor model.
\newblock {\em Rend. Sem. Mat. Univ. Politec. Torino} {\bf 57} (1999) 209--220.


\bibitem{SaVa}
{\sc F. Salvarani, J. L. V\' azquez.}
\newblock The diffusive limit for Carleman-type kinetic models.
\newblock {\em Nonlinearity} {\bf 18} (2005) 1223--1248.

\bibitem{Tartar}
{\sc L. Tartar.}
\newblock Some existence theorems for semilinear hyperbolic systems
in one space variable. 
\newblock Report \#2164, Mathematics Research Center,
University of Wisconsin, Madison, 1980.


\bibitem{Tay}
{\sc G. I. Taylor.}
\newblock Diffusion by continuous movements.
\newblock {\em Proc. London Math. Soc.} {\bf 20} (1922) 196--212.



\end{thebibliography}
\end{document}